\documentclass[a4paper, BCOR=0.0mm, DIV=calc]{scrartcl}
\usepackage[T1]{fontenc}
\usepackage[utf8]{inputenc}
\usepackage[ngerman]{babel}
\usepackage[osf,sc]{mathpazo}
\usepackage{textcomp}
\usepackage{microtype}
\usepackage{amsmath, amssymb, amsfonts}
\KOMAoptions{twoside=false, twocolumn=false, headinclude=false, footinclude=false, mpinclude=false, pagesize=auto}
\recalctypearea

\begin{document}
\title{\boldmath Über die Transformation der divergenten Reihe $1-mx+m(m+n)x^2-m(m+n)(m+2n)x^3$ + etc in einen Kettenbruch\footnote{
Originaltitel: "`De transformatione seriei divergentis $1 - mx + m(m+n)x^2 - m(m+n)(m+2n)x^3 + \mathrm{etc.}$ in fractionem continuam"', erstmals publiziert in "`\textit{Nova Acta Academiae Scientarum Imperialis Petropolitinae} 2, 1788, pp. 36-45 "', Nachdruck in "`\textit{Opera Omnia}: Series 1, Series 1, Volume 16, pp. 34 - 46 "', Eneström-Nummer E616, übersetzt von: Alexander Aycock, Textsatz: Artur Diener, im Rahmen des Projektes "`Eulerkreis Mainz"' }
\unboldmath}
\author{Leonhard Euler}
\date{}
\maketitle
\paragraph{§1}
Nachdem ich einst die Form divergenter Reihen solcher Art gründlicher untersucht hatte und den wahren Wert der hypergeometrischen Reihe
\[
	1 - 1 + 2 - 6 + 24 - 120 + 720 - \mathrm{etc.}
\]
mithilfe der Transformation in einen Kettenbruch angegeben hatte, kam ich auch auf diese sich viel weiter erstreckende Reihe
\[
	1 - mx + m(m+n)x^2 - m(m+n)(m+2n)x^3 + m(m+n)(m+2n)(m+3n)x^4 - \mathrm{etc.}
\]
zu sprechen, deren Summe ich gefunden hatte, diesem Kettenbruch
\[
	\cfrac{1}{1 + \cfrac{mx}{1 + \cfrac{nx}{1 + \cfrac{(m+n)x}{1 + \cfrac{2nx}{1+ \cfrac{(m+2n)x}{1+ \mathrm{etc.}}}}}}}
\]
gleich zu werden, die Gültigkeit welcher Sache ich aus der Verwandlung der Riccati-Gleichung in einen Kettenbruch gefolgert hatte. Weil aber dieser Beweis allzu weit hergeholt scheinen kann, werde ich hier dieselbe Reduktion aus einfacheren Prinzipien angeben.
\paragraph{§2}
Zuerst aber wird es nötig sein, diese allgemeine Reihe zu einer gefälligeren Form zusammenzuziehen, indem man
\[
	mx = a ~~ \mathrm{und} ~~ nx=b
\]
setzt, das diese unendliche Reihe vorgelegt ist
\[
	1 - a + a(a+b) - a(a+b)(a+2b) + a(a+b)(a+2b)(a+3b) - \mathrm{etc.}
\]
Damit aber außerdem die folgenden Auflösungen besser durchgeführt werden können und man nicht so viele Einschränkungen und Zusatzbedingungen braucht, möchte ich es setzen, wie folgt:
\[
	a = A,~ a+b=B,~ a+2b=C,~ a+3b=D,~ \mathrm{etc.}
\]
und so wird man diese Reihe haben
\[
	1 - A + AB - ABC + ABCD - \mathrm{etc.},
\]
deren gesuchte Summe wir mit dem Buchstaben $S$ bezeichnen wollen, so dass
\[
	S = 1 - A + AB - ABC + ABCD - \mathrm{etc.}
\]
ist; daher ist weiter
\[
	\frac{1}{S} = \frac{1}{1 - A + AB - ABC + ABCD - \mathrm{etc.}}
\]

\paragraph{§3}
Weil also $\frac{1}{S} > 1$ ist, führe man die letzte Gleichung zurück auf diese Form
\[
	\frac{1}{S} = 1 + \frac{A - AB + ABC - ABCD + \mathrm{etc.}}{1 - A + AB - ABC + ABCD - \mathrm{etc.}}.
\]
Nun aber wollen wir $\frac{1}{S} = 1 + \frac{1}{P}$ setzen und es wird
\[
	P = \frac{1 - A + AB - ABC + ABCD - \mathrm{etc.}}{1 - B + BC - BCD + BCDE - \mathrm{etc.}}
\]
sein; weil dieser Ausdruck wiederum die Einheit überträgt, und dieser wegen
\[
	B - A = b,~ C - A=2b,~ D-A=3b,~ \mathrm{etc.}
\]
sodann
\[
	P = 1 + \frac{b-2bB + 3bBC - 4bBCD + \mathrm{etc.}}{1 - B + BC - BCD + BCDE - \mathrm{etc.}}
\]
geben. Man setze also $P = 1 + \frac{b}{Q}$ und es wird
\[
	Q = \frac{1 - B + BC - BCD + BCDE - \mathrm{etc.}}{1 - 2B + 3BC - 4BCD + \mathrm{etc.}}
\]
sein, woher wir
\[
	Q = 1 + \frac{B - 2BC + 3BCD - 4BCDE + \mathrm{etc.}}{1 - 2B + 3BC - 4BCD + \mathrm{etc.}}
\]
folgern. Dieser Sache wegen wollen wir nun $Q = 1 + \frac{B}{R}$ setzen und es wird
\[
	R = \frac{1 - 2B + 3BC - 4BCD + \mathrm{etc.}}{1 - 2C + 3CD - 4CDE + \mathrm{etc.}}
\]
hervorgehen.
\paragraph{§4}
Hier tauchen also im Zähler wie im Nenner dieselben Koeffizienten auf, aber die Großbuchstaben sind im Nenner um einen Schritt vorgerückt worden. Weil also
\[
	C - B = b,~ D - B=2b,~ E - B = 3b,~ \mathrm{etc.}
\]
ist, wird
\[
	R = 1 + \frac{2b - 2\cdot 3bC + 3\cdot 4bCD - 4\cdot 5bCDE + \mathrm{etc.}}{1 - 2C + 3CD - 4CDE + 5CDEF - \mathrm{etc.}}
\]
werden. Wenn wir also $R = 1 + \frac{2b}{S}$ setzen, wird
\[
	S = \frac{1 - 2C + 3CD - 4CDE + \mathrm{etc.}}{1 - 3C + 6CD - 10CDE + \mathrm{etc.}}
\]
sein, wo im Nenner natürlich die Dreieckszahlen auftauchen; dieser Ausdruck wird zurückgeführt auf diesen
\[
	S = 1 + \frac{C - 3CD + 6CDE - 10CDEF + \mathrm{etc.}}{1 - 3C + 6CD - 10CDE + \mathrm{etc.}}.
\]
Wenn wir also $S = 1 + \frac{C}{T}$ setzen, wird
\[
	T = \frac{1 - 3C + 6CD - 10CDE + 15CDEF - \mathrm{etc.}}{1 - 3D + 6DE - 10DEF + 15DEFG - \mathrm{etc.}}
\]
sein.

\paragraph{§5}
Diese Form geht wegen
\[
	D-C=b,~ E-C=2b,~ F-C=3b,~ \mathrm{etc.}
\]
über in diese
\[
	T = 1 + \frac{3b-2\cdot 6bD + 3\cdot 10bDE - 4\cdot 15bDEF + \mathrm{etc.}}{1 - 3D + 6DE - 10DEF + 15DEFG - \mathrm{etc.}}
\]
Wir wollen $T = 1 + \frac{3b}{U}$ setzen, dass
\[
	U = \frac{1 - 3D + 6DE - 10DEF + 15DEFG - \mathrm{etc.}}{1 - 4D + 10DE - 20DEF + 35DEFG - \mathrm{etc.}}
\]
wird, wo im Nenner die ersten Pyramidalzahlen oder die Summen der Triagonalen gefunden werden, und daher erhalten wir
\[
	U = 1 + \frac{D - 4DE + 10DEF - 20DEFG + \mathrm{etc.}}{1 - 4D + 10DE - 20DEF + 35DEFG - \mathrm{etc.}},
\]
wo oben und unten schon die Pyramidalzahlen auftauchen. Man setze weiter
\[
	U = 1 + \frac{D}{V}
\]
und es wird
\[
	V = \frac{1 - 4D + 10DE - 20DEF + 35DEFG - \mathrm{etc.}}{1 - 4E + 10EF - 20EFG + 35EFGH - \mathrm{etc.}}
\]
werden.

\paragraph{§6}
Daher wird, indem man die Rechnung wie oben weiter verfolgt, weil
\[
	E - D = b,~ F-D=2b,~ G-D=3b,~ \mathrm{etc.}
\]
ist,
\[
	V = 1 + \frac{4b-2\cdot 10bE+3\cdot 20bEF - 4\cdot 35bEFG + \mathrm{etc.}}{1 - 4E + 10EF - 20EFG + 35EFGH - \mathrm{etc.}}
\]
sein. Es sei $V = \frac{4b}{X}$, dass 
\[
X = \frac{1- 4E + 10EF - 20EFG + 35EFGH - \mathrm{etc.} }{1 - 5E + 15EF - 35EFG + 70EFGH - \mathrm{etc.}}
\]
wird, welcher Ausdruck zurückgeführt auf diesen
\[
	X = 1 + \frac{E - 5EF + 15EFG - 35EFGH + \mathrm{etc.}}{1 - 5E + 15EF - 35EFG + \mathrm{etc.}}.
\]
Es sei $X = 1 + \frac{E}{Y}$ und es wird
\[
	Y = \frac{1 - 5E + 15EF - 35EFG + 70EFGH - \mathrm{etc.}}{1 - 5F + 15FG - 35FGH + 70FGHI - \mathrm{etc.}}
\]
sein.

\paragraph{§7}
Weil also
\[
	F - E = b,~ G - E = 2b,~ H - E = 3b,~ \mathrm{etc.}
\]
ist, wird
\[
	Y = 1 + \frac{5b - 2\cdot 15bF + 3\cdot 35bFG - 4\cdot 70bFGH + \mathrm{etc.}}{1 - 5F + 15FG - 35FGH + 70FGHI - \mathrm{etc.}}
\]
sein. Es sei nun $Y = 1 + \frac{5b}{Z}$, sodass
\[
	Z = \frac{1 - 5F + 15FG - 35FGH + 70FGHI - \mathrm{etc.}}{1 - 6F + 21FG - 56FGH + 126FGHI - \mathrm{etc.}}
\]
wird. Weil wir also eingangs $\frac{1}{S} = 1 + \frac{A}{P}$ gesetzt haben, wird die gesuchte Summe
\[
	S = \cfrac{1}{1 + \cfrac{A}{P}}
\]
sein; dann aber sind die folgenden Festsetzungen gemacht worden:
\begin{align*}
	 P = 1 + \frac{b}{Q},~~~~ S = 1 + \frac{C}{T},~~~~ V = 1 + \frac{4b}{X} \\
	 Q = 1 + \frac{B}{R},~~~~ T = 1 + \frac{3b}{U},~~~ X = 1 + \frac{E}{Y} \\
	 R = 1 + \frac{2b}{S},~~~ U = 1 + \frac{D}{V},~~~~ Y = 1 + \frac{5b}{Z}\\
	\mathrm{etc.}
\end{align*}
Nachdem der Reihe diese Werte eingesetzt worden sind, entsteht dieser Kettenbruch:
\[
	S = \cfrac{1}{1 + \cfrac{A}{1 + \cfrac{b}{1 + \cfrac{B}{1 + \cfrac{2b}{1 + \cfrac{C}{1 + \cfrac{3b}{1 + \cfrac{D}{1 + \cfrac{4b}{1 + \mathrm{etc.}}}}}}}}}}
\]
Wenn wir also anstelle der Buchstaben $A,\, B,\, C,\, D$ wieder die angenommenen Werte einsetzen, dass wir diese divergente Reihe haben
\[
	1 - a + a(a+b) - a(a+b)(a+2b) + a(a+b)(a+2b)(a+3b) - \mathrm{etc},
\]
deren Summe durch den folgenden Kettenbruch ausgedrückt werden wird:
\[
	S = \cfrac{1}{1 + \cfrac{a}{1 + \cfrac{b}{1 + \cfrac{a+b}{1 + \cfrac{2b}{1 + \cfrac{a+2b}{1 + \cfrac{3b}{1 + \cfrac{a+3b}{1 + \cfrac{4b}{1 + \mathrm{etc.}}}}}}}}}}
\]
welches dieselbe Form ist, die ich einst gegeben hatte.

\paragraph{§8}
Diese Tranformation ist umso bemerkenswerter, weil sie uns einen sehr sicheren und den vielleicht einzigen Weg eröffnet, den Wert einer divergenten Reihe zumindest näherungsweise zu bestimmen. Wenn nämlich ein Kettenbruch auf gewohnte Weise in einfachere Brüche aufgelöst wird,
\[
	1,~ \frac{1}{1+a},~ \frac{1+b}{1+a+b},~ \mathrm{etc,}
\]
sind diese abwechselnd größer und kleiner als der Wert der divergenten Reihe und gehen immer näher an diesen Wert heran. Dann aber habe ich auch einst einzigartige Kunstgriffe ausgelegt, die um Vieles schneller zum wahren Wert führen.

\paragraph{§9}
Außerdem aber wird es auch förderlich sein, einen solchen Kettenbruch
\[
	1 + \cfrac{\alpha}{1 + \cfrac{\beta}{1 + \cfrac{\gamma}{1 + \cfrac{\delta}{1 + \mathrm{etc.}}}}}
\]
im Allgemeinen hinreichend angenehm auf die Anzahl der Teile zurückführen zu können. Nachdem nämlich sein Wert gleich $S$ gesetzt worden ist, wird er sich so darstellen lassen:
\[
	S = 1 + \cfrac{\alpha}{1 + \cfrac{\beta}{P}},~~ P = 1 + \cfrac{\gamma}{1 + \cfrac{\delta}{Q}},~~ Q = 1 + \cfrac{\varepsilon}{1 + \cfrac{\zeta}{R}},~~ \mathrm{etc.}
\]
Schon die erste dieser Formeln wird
\[
	S = 1 + \cfrac{\alpha P}{P + \beta} = 1 + \alpha - \cfrac{\alpha \beta}{\beta + P}
\]
sein, die zweite Formel gibt darauf
\[
	P = 1 + \cfrac{\gamma Q}{Q + \delta} = 1 + \gamma - \cfrac{\gamma \delta}{\delta + Q},
\]
auf dieselbe Weise liefert die dritte
\[
	Q = 1 + \cfrac{\varepsilon R}{R + \zeta} = 1 + \varepsilon - \cfrac{\varepsilon \zeta}{\zeta + R}, \mathrm{etc.}
\]
Diese nacheinander eingesetzten Werte werden daher diesen neuen Kettenbruch ergeben:
\[
	S = 1 + \alpha - \cfrac{\alpha \beta}{1 + \beta + \gamma - \cfrac{\gamma\delta}{1+\delta + \varepsilon - \cfrac{\varepsilon\zeta}{1+\zeta +\eta - \cfrac{\eta \Theta}{1 + \Theta + i - \mathrm{etc.}}}}}
\]

\paragraph{§10}
Weil daher in unserem Fall die divergente Reihe
\[
	S = 1 - a + a(a+b) - a(a+b)(a+2b) + a(a+b)(a+2b)(a+3b) - \mathrm{etc.}
\]
auf diesen Kettenbruch
\[
	S = \cfrac{1}{1 + \cfrac{a}{1 + \cfrac{b}{1 + \cfrac{a+b}{1 + \cfrac{2b}{1 + \cfrac{a+2b}{1 + \cfrac{3b}{1 + \cfrac{a+3b}{1 + \mathrm{etc.}}}}}}}}}
\]
zurückgeführt worden ist, wollen wir hier $\alpha = a,\, \beta = b,\, \gamma = a+b,\, \delta = 2b,\, \varepsilon = a+2b,\, \mathrm{etc.}$ nehmen und es wird
\[
	S = 1 + a - \cfrac{ab}{1+a+2b - \cfrac{2b(a+b)}{1+a+4b - \cfrac{3b(a+2b)}{1+a+6b - \cfrac{4b(a+3b)}{1+a+\mathrm{etc.}}}}}
\]
sein.

\appendix
\section*{Über den Brouncker'schen Kettenbruch}
\paragraph{§11}
Nachdem ich einst viel beschäftigt war mit dem Erforschen der Analysis, die Brouncker zu diesem einzigartigen Kettenbruch geführt hat, weil es mir natürlich unwahrscheinlich schien, dass er durch so viele Umwege, wie von Wallis erwähnt werden, dorthin geführt worden ist, schien es schließlich für mich, schon hinreichend klar gezeigt zu haben, dass Brouncker diese Form aus der Leibniz'schen Reihe
\[
	1 - \frac{1}{3} + \frac{1}{5} - \frac{1}{7} + \frac{1}{9} - \frac{1}{11} + \mathrm{etc,}
\]
welche der große Gregory schon vorher gefunden hatte, gefolgert hat, eher als aus der Interpolation der Reihe
\[
	1,~~ \frac{1}{2},~~ \frac{1\cdot 3}{2\cdot 4},~~ \frac{1\cdot 3 \cdot 5}{2\cdot 4 \cdot 6},~~ \frac{1\cdot 3\cdot 5 \cdot 7}{2\cdot 4 \cdot 6\cdot 8},~~ \mathrm{etc,}
\]
so wie es Wallis vermutete, weil natürlich die Betrachtung jener Reihe durch hinreichend klare Rechnung auf die Brouncker'sche Form führt.

\paragraph{§12}
Diese Beobachtung scheint nun aber freilich umso größerer Aufmerksamkeit würdig, nachdem der berühmte Daniel Bernoulli die Erwähnung der Brouncker'schen Form zu erneuern nicht unwürdig gefunden hatte. Weil ich ja also vor nicht so langer Zeit eine leichte Methode erörtert habe, diese Form aus der Reihe
\[
	1 - \frac{1}{3} + \frac{1}{5} - \frac{1}{7} + \mathrm{etc.}
\]
zu berechnen, glaube ich, dass es für die Mathematiker nicht unangemessen sein wird, wenn ich die inverse Methode in den Mittelpunkt stellen werde, mit deren Hilfe sich umgekehrt die Brouncker'sche Form auf die Leibniz'sche Reihe zurückführen lässt.

\paragraph{§13}
Ich werde also diesen Kettenbruch betrachten, als ob sein Wert noch nicht bekannt wäre, indem ich
\[
	S = \cfrac{1}{1 + \cfrac{1}{2 + \cfrac{9}{2 + \cfrac{25}{2 + \cfrac{49}{2 + \cfrac{81}{2 + \mathrm{etc.}}}}}}}
\]
setze, welchen ich durch Teile auf folgende Weise darstelle:
\[
	S = \cfrac{1}{1 + \cfrac{1}{-1+P}},~~ P = 3 + \cfrac{9}{-3 + Q},~~ Q = 5 + \cfrac{25}{-5 + R},~~ R = 7 + \cfrac{49}{-7 + S},~~ \mathrm{etc.}
\]
Aus diesen entsprechend verbundenen Teilen nämlich entsteht natürlich die vorgelegte Form selbst.

\paragraph{§14}
Wir wollen also diese Teile getrennt entwickeln; und freilich liefert die erste auf einen einfachen Bruch zurückgeführt
\[
	S = \frac{P-1}{P}
\]
und daher $S = 1- \frac{1}{P}$; die zweite aber wird $\frac{3Q}{Q-3}$ sein, woher
\[
	\frac{1}{P} = \frac{Q - 3}{3Q} ~~ \mathrm{oder} ~~ \frac{1}{P} = \frac{1}{3} - \frac{1}{Q}
\]
sein wird; auf ähnliche Weise gibt der dritte Teil
\[
	Q = \frac{5R}{R-5}
\]
und daher $\frac{1}{Q} = \frac{1}{5} - \frac{1}{R}$; auf dieselbe Weise werden wir aus den folgenden Teilen
\[
	\frac{1}{R} = \frac{1}{7} - \frac{1}{5},~~ \frac{1}{S} = \frac{1}{9} - \frac{1}{T},~ \mathrm{etc.}
\]
erhalten. Wenn daher diese Werte nacheinander eingesetzt werden, werden wir diesen Ausdruck erhalten
\[
	S = 1 - \frac{1}{3} + \frac{1}{5} - \frac{1}{7} + \frac{1}{9} - \frac{1}{11} + \frac{1}{13} - \mathrm{etc,}
\]
so dass wir nun nun gewiss sind, dass $S = \frac{\pi}{4}$ ist.

\paragraph{§15}
Auf ähnliche Weise wird sich auch der Wert anderer Kettenbrüche solcher Art untersuchen lassen. Wenn z.\,B. diese Form vorgelegt war
\[
	S = \cfrac{1}{1 + \cfrac{1}{1 + \cfrac{4}{1 + \cfrac{9}{1 + \cfrac{16}{1 + \mathrm{etc.}}}}}}
\]
teile man diese auf die folgende Weise in Gliedern auf:
\[
	S = \cfrac{1}{1 + \cfrac{1}{-1+P}},~~ P = 2 + \cfrac{4}{-2 + Q},~~ Q = 3 + \cfrac{9}{-3 + R},~~ R = 4 + \cfrac{16}{-4 + S},~~ \mathrm{etc;}
\]
nachdem nämlich diese einzelnen Teile entwickelt worden sind, wird man
\[
	S = 1 - \frac{1}{P},~~ \frac{1}{P} = \frac{1}{2} - \frac{1}{Q},~~ \frac{1}{Q} = \frac{1}{3} - \frac{1}{R},~~ \frac{1}{R} = \frac{1}{4} - \frac{1}{S},~~ \mathrm{etc.}
\]
finden, woher man folgt, dass
\[
	S = 1 - \frac{1}{2} + \frac{1}{3} - \frac{1}{4} + \frac{1}{5} - \frac{1}{6} + \mathrm{etc.} = \log{2}
\]
sein wird. Diese Methode also scheint genauso beim Zurückgehen zu gelten.
\end{document}